\documentclass[twoside,11pt]{amsart}
\usepackage{latexsym}
\usepackage{graphics, epsfig}
\newtheorem{theorem}{Theorem}[section]

\newtheorem{corollary}{Corollary}[section]

\begin{document}

\title{Exotic smooth structures on
       $3{\mathbf CP}^2 \sharp 8{\overline{{\mathbf CP}}^2}$}

\author{Jongil Park*}

\address{Department of Mathematical Sciences, Seoul National University\\
         San 56-1 Sillim-dong, Gwanak-gu, Seoul 151-747, Korea}

\email{jipark@math.snu.ac.kr}

\thanks{* This work was supported by grant No. R01-2005-000-10625-0 from
        the KOSEF}

\thanks{* The author also holds a joint appointment in the Research Institute
        of Mathematics, SNU}




\begin{abstract}
 Motivated by Stipsicz and Szab\'{o}'s exotic $4$-manifolds with
 $b_2^+=3$ and $b_2^-=8$, we construct a family of
 simply connected smooth $4$-manifolds with $b_2^+=3$ and $b_2^-=8$.
 As a corollary, we conclude that the topological $4$-manifold
 $3{\mathbf CP}^2 \sharp 8{\overline{{\mathbf CP}}^2}$ admits
 infinitely many distinct smooth structures.
\end{abstract}

\maketitle

\section{Introduction}

\markboth{JONGIL PARK}{EXOTIC SMOOTH STRUCTURES ON
       $3{\mathbf CP}^2 \sharp 8{\overline{{\mathbf CP}}^2}$}

 Since S. Donaldson introduced gauge theory in 1982, it has been known
 that most simply connected smooth $4$-manifolds with $b_2^+$ odd
 and large enough admit infinitely many distinct smooth structures
 (\cite{fs2}, \cite{g}, \cite{p2}, \cite{p3}).
 But it is still an intriguing problem to determine which smooth
 $4$-manifolds with $b_2^+$ small have more than one smooth structure.
 In the case when $b_2^+ = 1$, S. Donaldson first proved that
 a Dolgachev surface is not diffeomorphic to
 ${\mathbf CP}^2 \sharp 9{\overline{{\mathbf CP}}^2}$ (\cite{d})
 and D. Kotschick proved in the late 1980's that
 the Barlow surface is not diffeomorphic to
 ${\mathbf CP}^2 \sharp 8{\overline{{\mathbf CP}}^2}$ (\cite{k}).
 Recently, the author constructed a new simply connected symplectic
 $4$-manifold with $b_2^+=1$ and $b_2^-=7$ (\cite{p4}),
 and then R. Fintushel, R. Stern, A. Stipsicz
 and Z. Szab\'{o} found many new exotic
 smooth $4$-manifolds with $b_2^+=1$ using a rational
 blow-down technique (\cite{fs3}, \cite{pss}, \cite{ss2}).
 So it has been proved up to now that rational surfaces
 ${\mathbf CP}^2 \sharp n{\overline{{\mathbf CP}}^2}$ with $n \geq 5$
 admit infinitely many distinct smooth structures.
 On the other hand, the case with $b_2^+=3$ has also been studied extensively
 for a long time. For example, it was known in the mid 1990's
 that the $K3$ surface
 $E(2)$ and $3{\mathbf CP}^2 \sharp 19{\overline{{\mathbf CP}}^2}$
 admit infinitely many distinct smooth structures (\cite{fs1},
 \cite{fm}, \cite{mo}, \cite{ss1}) and the topological $4$-manifolds
 $3{\mathbf CP}^2 \sharp n{\overline{{\mathbf CP}}^2}$ with $n \geq 14$
 also admit infinitely many distinct smooth structures (\cite{st}, \cite{sz}).
 And later, the same statement with $n \geq 10$ was also proved
 (\cite{pd1}, \cite{pd2}, \cite{pd3}). Recently,
 Stipsicz and Szab\'{o} constructed infinitely
 many distinct smooth structures on
 $3{\mathbf CP}^2 \sharp 9{\overline{{\mathbf CP}}^2}$
 and they also constructed infinitely many exotic smooth $4$-manifolds
 with vanishing first homology, $b_2^+ =3$ and $b_2^- =8$ using a similar
 technique appeared in~\cite{pss} (\cite{ss3}).

 The aim of this paper is to prove that the topological $4$-manifold
 $3{\mathbf CP}^2 \sharp 8{\overline{{\mathbf CP}}^2}$ also admits
 infinitely many distinct smooth structures.
 We first construct a family of simply connected smooth
 $4$-manifolds by modifying Stipsicz and Szab\'{o}'s exotic
 $4$-manifolds with $b_2^+=3$ and $b_2^-=8$ in~\cite{ss3}.
 And then we compute their Seiberg-Witten invariants to show that
 they are mutually nondiffeomorphic.

\begin{theorem}
\label{thm-main}
 The simply connected topological $4$-manifold
 $3{\mathbf CP}^2 \sharp 8{\overline{{\mathbf CP}}^2}$
 admits infinitely many distinct smooth structures.
\end{theorem}

\noindent {\bf Acknowledgement}.
 The author would like to thank Andr\'{a}s  Stipsicz for his
 invaluable help and advice while working on this problem.
 He also pointed out critical errors in the early versions
 of this paper to the author.

\bigskip

\section{A Main Construction}
\label{sec-2}

 In this section we first briefly review a rational blow-down surgery
 and state related theorems (\cite{fs1}, \cite{p1} for details).
 Suppose that $p>q>0$ are relatively prime integers.
 Let $C_{p,q}$ be a smooth $4$-manifold obtained by plumbing disk
 bundles over the $2$-sphere instructed by the following linear diagram

\begin{picture}(400,40)(-100,-20)
  \put(0,3){\makebox(200,20)[bl]{$-b_{k}$ \hspace{15pt} $-b_{k-1}$
                                  \ \ \hspace{70pt}  $-b_{1}$}}
  \put(4,-25){\makebox(200,20)[tl]{ $u_{k}$ \hspace{21pt} $u_{k-1}$
        \hspace{85pt} $u_{1}$}}  \multiput(10,0)(40,0){2}{\line(1,0){40}}
                         \multiput(10,0)(40,0){2}{\circle*{3}}
                         \multiput(100,0)(5,0){4}{\makebox(0,0){$\cdots$}}
  \put(125,0){\line(1,0){40}}
  \put(165,0){\circle*{3}}
\end{picture}

\noindent
 where $\frac{p^{2}}{pq-1} =[b_{k},b_{k-1}, \ldots, b_{1}]$
 is the unique continued fraction with all $b_{i} \geq 2$
 and each vertex represents a disk bundle over the $2$-sphere
 $u_{i}$ whose Euler number is $-b_{i}$.
 Then $C_{p,q}$ is a negative definite simply connected $4$-manifold
 whose boundary is a lens space $L(p^2, 1-pq)$.
 It was known that a lens space $L(p^2, 1-pq)$ bounds a rational ball
 $B_{p,q}$ with $\pi_{1}(B_{p,q}) \cong {\mathbf Z}_{p}$ (\cite{ch}).
 Furthermore, the inclusion
 $\partial B_{p,q} \longrightarrow B_{p,q}$
 induces an epimorphism
 $\pi_{1}(\partial B_{p,q})\cong {\mathbf Z}_{p^2}
 \longrightarrow \pi_{1}(B_{p,q}) \cong {\mathbf Z}_{p}$. \\

 \noindent {\bf Definition} Suppose $X$ is a smooth $4$-manifold
 containing a configuration $C_{p,q}$. Then we may construct a new smooth
 $4$-manifold $X_{p,q}$, called the {\em (generalized) rational blow-down} of $X$,
 by replacing $C_{p,q}$ with the rational ball $B_{p,q}$.
 Note that this process is well-defined, that is, a new smooth
 $4$-manifold $X_{p,q}$ is uniquely determined (up to diffeomorphism)
 from $X$ because each diffeomorphism of $\partial B_{p,q}$
 extends over the rational ball $B_{p,q}$.
 We call this a {\em rational blow-down} surgery. \\

 A rational blow-down surgery with $q=1$ was originally
 introduced by Fintushel and Stern to compute the Donaldson series
 for simply connected regular elliptic surfaces with multiple fibers of
 relatively prime orders (\cite{fs1}) and later it was extended to
 the general case  by the author (\cite{p1}). Furthermore,
 it was proved that the Seiberg-Witten basic classes and the Seiberg-Witten
 invariants of a smooth $4$-manifold $X_{p,q}$ obtained by a rational
 blow-down surgery are closely related to those of $X$.
 For example, we have the following results.

\begin{theorem}[\cite{p1}]
\label{thm-sw}
 Suppose $X$ is a smooth $4$-manifold with $b_2^+ >1$
 which contains a configuration $C_{p,q}$.
 If $L$ is a SW-basic class of $X$ which satisfies
 $(L|_{C_{p,q}})^{2} =-b_{2}(C_{p,q})$ and
 $c_{1}(L|_{L(p^2, 1-pq)})  = mp \in {\mathbf Z}_{p^{2}} \cong
 H^{2}(L(p^2, 1-pq);{\mathbf Z})$ with $m \equiv (p-1)\!\pmod2$,
 then $L$ induces a SW-basic class $\overline{L}$ of $X_{p,q}$
 such that $SW_{X_{p,q}}(\overline{L}) = SW_{X}(L)$.
\end{theorem}

\begin{corollary}
\label{cor-sw}
 Suppose $X$ is a smooth $4$-manifold with $b_2^+ >1$
 which contains a configuration $C_{p,q}$.
 If $L$ is a SW-basic class of $X$ satisfying
 $L\cdot u_{i} = (b_{i}-2),\, \forall \, 1 \leq i \leq k$
 $($or $L\cdot u_{i} =-(b_{i}-2),\, \forall\, i )$,
 then $L$ induces a SW-basic class $\overline{L}$ of $X_{p,q}$
 such that $SW_{X_{p,q}}(\overline{L}) = SW_{X}(L)$.
\\
\end{corollary}

 Next, we review the exotic smooth $4$-manifolds with $b_2^+=3$ and
 $b_2^-=8$ constructed by Stipsicz and Szab\'{o}
 (\cite{ss3} for details).

 Let $E(2)$ be the simply connected elliptic surface with
 no multiple fibers and holomorphic Euler characteristic $2$.
 Then $E(2)$ is an elliptic fibration over $S^2$ which has
 various types of singular fibers.
 In particular, it can be viewed as an elliptic fibration over $S^2$
 with an $I_{16}$-singular fiber and eight fishtail-singular fibers
 (\cite{hkk}). Note that one can change the eight fishtail fibers
 to three pairs of the same type fishtail fibers and
 two other fishtail fibers in $E(2)$.
%
%
 Let $V_{K_1, K_2, K_3}$ be a smooth $4$-manifold obtained by doing
 three knot surgeries in the three double node neighborhoods,
 each of which contains a pair of the same type fishtail fibers,
 in $E(2)$ respectively, using the twist knots $K_1$ for the first,
 $K_2$ for the second and $K_3$ for the third surgery.
 Then
 the smooth $4$-manifold $V_{K_1, K_2, K_3}$ contains
 a `pseudo-section' i.e. an immersed $2$-sphere with three positive
 double points and with square $-2$ (\cite{fs3}).
 Furthermore,
 resolving the two positive intersections of this pseudo-section
 with the remaining two fishtail fibers,
 the $4$-manifold $V_{K_1, K_2, K_3}$ contains an immersed $2$-sphere
 with five positive double points and with square $2$.
 By blowing up at the five positive double points of the
 immersed $2$-sphere,
 $V_{K_1, K_2, K_3}\sharp 5{\overline{{\mathbf CP}}^2}$
 contains an embedded $2$-sphere $S$ of square $-18$ which
 still intersects the $I_{16}$-fiber transversally
 at one point. Let $\Sigma_{0}$ be an embedded $2$-sphere
 in $I_{16}$-fiber which intersects $S$ and let $\Sigma_{1}$
 be an embedded $2$-sphere in $I_{16}$-fiber which intersects with
 $\Sigma_{0}$ at one point, say $q$, positively. Then, applying $17$
 infinitely close blow-ups at the point $q$,
 the $4$-manifold $V_{K_1, K_2, K_3}\sharp 22{\overline{{\mathbf CP}}^2}$
 contains a configuration $C_{305, 17}$ which is a chain of embedded
 $2$-spheres according to the following linear plumbing
 \[ (-18, -19, \, \overbrace{-2,\, \ldots, -2}^{14}, \, -3, \,
   \overbrace{-2,\, \ldots, -2}^{16}\, ). \]

 Note that the boundary $\partial{C_{305, 17}}$ is the lens space
 $L(305^2, -5184)$ which also bounds a rational ball $B_{305, 17}$.
 Hence,
 by rationally blowing down along the configuration $C_{305, 17}$
 in $V_{K_1, K_2, K_3}\sharp 22{\overline{{\mathbf CP}}^2}$,
 Stipsicz and Szab\'{o} constructed a new family of smooth $4$-manifolds
 \[ Y_{K_1, K_2, K_3}= (V_{K_1, K_2, K_3}\sharp 22{\overline{{\mathbf
 CP}}^2} - \mathrm{int}\, C_{305, 17}) \cup_{L(305^2,-5184)} B_{305, 17}.
 \]
 Furthermore,
 they computed the Seiberg-Witten invariants of these manifolds.
 In particular, they showed that an infinite family of smooth
 $4$-manifolds $Y_{n}:= Y_{T_n, T_n, T_n}$ are pairwise nondiffeomorphic,
 all with nonvanishing Seiberg-Witten invariants.

\begin{theorem}[\cite{ss3}]
\label{thm-ss3}
 There are infinitely many pairwise nondiffeomorphic smooth, closed
 $4$-manifolds $Y_n$ with vanishing first homology, $b_2^+=3, \,
 b_2^-=8$ and nontrivial Seiberg-Witten invariants.
\end{theorem}

 Now we will modify their construction to get an infinite family
 of simply connected smooth $4$-manifolds with $b_2^+=3$ and
 $b_2^-=8$. The procedure is the following:
 Recall that the $4$-manifold $V_{K_1, K_2, K_3}$ above
 contains a `pseudo-section' i.e. an immersed $2$-sphere
 with three positive double points and with square $-2$.
 At this stage,
 instead of resolving two positive intersections of this
 pseudo-section with the remaining two fishtail fibers,
 resolve just one positive intersection point.
 Then the $4$-manifold $V_{K_1, K_2, K_3}$ contains an immersed
 $2$-sphere with four positive double points and with square $0$.
 By blowing up at the four positive double points of the immersed
 $2$-sphere,
 $V_{K_1, K_2, K_3}\sharp 4{\overline{{\mathbf CP}}^2}$
 contains an embedded $2$-sphere of square $-16$ which still
 intersects one remaining fishtail fiber at one point, say $p$,
 and also intersects transversally the $I_{16}$-fiber.
 Now applying $2$ infinitely close blow-ups at the point $p$ and
 blowing up at the positive double point lying in the remaining
 fishtail fiber, we get a $4$-manifold
 $V_{K_1, K_2, K_3}\sharp 7{\overline{{\mathbf CP}}^2}$
 which contains an embedded $2$-sphere $S$ of square $-18$
 and a configuration $C_{3,1}$ disjoint from $S$.
 And then, proceeding as above, i.e. applying $17$
 infinitely close blow-ups at the point $q$,
 we obtain a smooth $4$-manifold
 $V_{K_1, K_2, K_3}\sharp 24{\overline{{\mathbf CP}}^2}$
 which contains two disjoint configurations $C_{3,1}$ and $C_{305, 17}$.
 More precisely,
 the $4$-manifold $V_{K_1, K_2, K_3}\sharp 24{\overline{{\mathbf CP}}^2}$
 contains a chain of embedded $2$-spheres according to the following
 linear plumbing
\begin{eqnarray}
 (-5, -2, -1, -18, -19, \, \overbrace{-2,\, \ldots, -2}^{14}, \, -3, \,
   \overbrace{-2,\, \ldots, -2}^{16}\, ).
\end{eqnarray}

 Note that the third embedded $2$-sphere, say $E_{6}$, of square $-1$
 in the linear plumbing above is an exceptional curve coming from
 the $6^{th}$ blowing up of $V_{K_1, K_2, K_3}$.
 Since two configurations $C_{3, 1}$ and $C_{305, 17}$ are disjoint
 in $V_{K_1, K_2, K_3}\sharp 24{\overline{{\mathbf CP}}^2}$,
 we first get a family of smooth $4$-manifolds, say $Y'_{K_1, K_2, K_3}$,
 by rationally blowing down along $C_{305, 17}$,
 and then we finally get an infinite family of new smooth $4$-manifolds,
 say $Z_{K_1, K_2, K_3}$, by rationally blowing down along $C_{3, 1}$.
 Alternatively, we can get a family of smooth $4$-manifolds,
 say $Y''_{K_1, K_2, K_3}$, by rationally blowing down along $C_{3, 1}$
 first, and then we can get the same smooth $4$-manifolds $Z_{K_1, K_2, K_3}$
 by rationally blowing down along $C_{305, 17}$. That is, we have
\begin{eqnarray}
 Y'_{K_1, K_2, K_3} &=& (V_{K_1, K_2, K_3}\sharp 24{\overline{{\mathbf
    CP}}^2} - \mathrm{int}\, C_{305, 17}) \cup_{L(305^2,-5184)} B_{305, 17} \\
 Y''_{K_1, K_2, K_3} &=& B_{3, 1} \cup_{L(9,-2)}
    (V_{K_1, K_2, K_3}\sharp 24{\overline{{\mathbf CP}}^2} - \mathrm{int}\, C_{3, 1}) \\
 Z_{K_1, K_2, K_3} &=& B_{3, 1} \cup_{L(9,-2)}
    (Y'_{K_1, K_2, K_3} - \mathrm{int}\, C_{3, 1}) \\
    &=& (Y''_{K_1, K_2, K_3} - \mathrm{int}\, C_{305, 17})
    \cup_{L(305^2, -5184)} B_{305, 17}.
\end{eqnarray}

 Next,
 we will prove that the $4$-manifolds $Z_{K_1, K_2, K_3}$ constructed
 above are simply connected by using Van-Kampen's theorem.

\begin{theorem}
\label{thm-main-sc}
 For any twisted knots $K_1, K_2$ and $K_3$, the $4$-manifold
 $Z_{K_1, K_2, K_3}$ is homeomorphic to
 $3{\mathbf CP}^2 \sharp 8{\overline{{\mathbf CP}}^2}$.
\end{theorem}

\begin{proof}
 First, in order to prove the simple connectivity of
 $Z_{K_1, K_2, K_3}$, let us decompose the smooth $4$-manifold
 $V_{K_1, K_2, K_3}\sharp 24{\overline{{\mathbf CP}}^2}$ into
 $C_{3,1}\cup_{L(9,-2)} X_{0} \cup_{L(305^2, -5184)} C_{305,17}$.
 Then the rational blow-down $4$-manifolds $Y'_{K_1, K_2, K_3}$,
 $Y''_{K_1, K_2, K_3}$ and $Z_{K_1, K_2, K_3}$
 can be decomposed as follows:
\begin{eqnarray}
 Y'_{K_1, K_2, K_3} &=& C_{3,1}\cup_{L(9,-2)} X_{0} \cup_{L(305^2, -5184)} B_{305,17} \\
 Y''_{K_1, K_2, K_3} &=& B_{3,1}\cup_{L(9,-2)} X_{0} \cup_{L(305^2, -5184)} C_{305,17} \\
 Z_{K_1, K_2, K_3} &=& B_{3,1}\cup_{L(9,-2)} X_{0} \cup_{L(305^2, -5184)} B_{305,17}.
\end{eqnarray}

\noindent
 Let $i'_{*}: \pi_{1}(\partial C_{305,17}) \rightarrow
 \pi_{1}(C_{3,1} \cup_{L(9,-2)}X_{0})$ and
 $i''_{*}: \pi_{1}(\partial C_{3,1}) \rightarrow
 \pi_{1}(X_{0} \cup_{L(305^2, -5184)} C_{305,17})$
 be induced homomorphisms by inclusions
 $i': \partial C_{305,17} \rightarrow C_{3,1}\cup_{L(9,-2)} X_{0}$
 and
 $i'': \partial C_{3,1} \rightarrow X_{0} \cup_{L(305^2, -5184)} C_{305,17}$
 respectively.
 We may also choose the generators, say $\alpha$ and $\beta$, of
 $\pi_{1}(\partial C_{305,17}, x_{0}) \cong {\mathbf Z}_{305^2}$
 and $\pi_{1}(\partial C_{3,1}, x_{1}) \cong {\mathbf Z}_{9}$
 so that $\alpha$ and $\beta$ are represented by circles lying in
 $\partial C_{305,17} \cap E_6$ and $\partial C_{3,1} \cap E_6$,
 respectively.
 Recall that $E_6$ is an exceptional curve which intersects
 transversally each embedded $2$-sphere ending in
 the configurations $C_{3,1}$ and $C_{305,17}$
 (refer to the linear plumbing $(1)$ above).
 Since $V_{K_1, K_2, K_3}\sharp 24{\overline{{\mathbf CP}}^2}$
 is simply connected and since the images $i'_{*}(\alpha)$ and
 $i''_{*}(\beta)$ of generators are bounded by disks coming from
 a punctured exceptional curve
 $E_6\setminus \{\mathrm{an\, \, open\, \, disk}\}$,
 the two $4$-manifolds $Y'_{K_1, K_2, K_3}$ and $Y''_{K_1, K_2, K_3}$
 are simply connected.
 Hence, applying Van-Kampen's theorem on the decompositions
 $(6)$ and $(7)$ above,
 we conclude that both quotients groups
\begin{eqnarray}
  \pi_{1}(X_{0} \cup_{L(305^2,-5184)} B_{305,17}, x_{1})/N_{j'_{*}(\beta)}\,
    \, \, \, \mathrm{and}\, \, \, \,
    \pi_{1}(B_{3,1} \cup_{L(9, -2)} X_{0}, x_{0})/N_{j''_{*}(\alpha)}
\end{eqnarray}

\noindent
 are trivial. Here
 $j': \partial C_{3,1} \rightarrow X_{0} \cup_{L(305^2,-5184)} B_{305,17}$
 and
 $j'': \partial C_{305,17} \rightarrow B_{3,1} \cup_{L(9, -2)} X_{0}$
 are inclusions, and
 $N_{j'_{*}(\beta)}$ and $N_{j''_{*}(\alpha)}$ are the least
 normal subgroups containing $j'_{*}(\beta)$ and $j''_{*}(\alpha)$
 respectively.
 Note that they satisfy ${j'_{*}(\beta)}^9 = 1$ and
 ${j''_{*}(\alpha)}^{305^2} = 1$. Furthermore, there is a relation
 between $j'_{*}(\beta)$ and $j''_{*}(\alpha)$ when we restrict them to
 $X_{0}$. That is, if we choose a path $\gamma$ connecting $x_0$ and $x_1$
 lying in $E_6 \setminus \{\mathrm{two\, \, open\, \, disks}\}$,
 then they satisfy either
 $j'_{*}(\beta) = \gamma^{-1} \cdot j''_{*}(\alpha)\cdot \gamma$
 or  $j'_{*}(\beta) = \gamma^{-1} \cdot j''_{*}(\alpha)^{-1} \cdot \gamma$
 (depending on orientations) because one is homotopic to the other in
 $E_6 \setminus \{\mathrm{two\, \, open\, \, disks} \} \subset X_{0}$.
 Hence, by combining two facts above, we get
 ${j'_{*}(\beta)}^{305^2} = (\gamma^{-1} \cdot j''_{*}(\alpha)^{\pm 1} \cdot
 \gamma)^{305^2} = \gamma^{-1} \cdot j''_{*}(\alpha)^{\pm 305^2} \cdot \gamma
  = 1 = {j'_{*}(\beta)}^9$. Since the two numbers $9$ and $305^2 = 25 \times
  61^2$ are relatively prime, the element ${j'_{*}(\beta)}$ should
  be trivial. So the relation
  $j'_{*}(\beta) = \gamma^{-1} \cdot j''_{*}(\alpha)^{\pm 1} \cdot \gamma$
  implies the triviality of $j''_{*}(\alpha)$.
  Finally the relation $(9)$ above implies that both groups
  \begin{eqnarray}
  \pi_{1}(X_{0} \cup_{L(305^2,-5184)} B_{305,17})\,
    \, \, \, \mathrm{and}\, \, \, \,
    \pi_{1}(B_{3,1} \cup_{L(9, -2)} X_{0})
\end{eqnarray}

\noindent
 are trivial. Thus, applying Van-Kampen's theorem again on
 the decomposition $(8)$ above, we conclude that the $4$-manifold
 $Z_{K_1, K_2, K_3}$ is simply connected.

 The rest of proof follows from simple Euler characteristics and
 signature computations together with Freedman's classification
 theorem.
\end{proof}

 Finally, we compute the Seiberg-Witten invariants of smooth $4$-manifolds
 $Z_{K_1, K_2, K_3}$ to show that they are mutually nondiffeomorphic.
 For simplicity, we restrict our attention to the special case
 $K_1 = K_2 = K_3 = n$-twist knot $T_n$.
 Let $V_n, Y'_{n}, Y''_{n}$ and $Z_{n}$ denote simply connected smooth
 $4$-manifolds $V_{T_n, T_n, T_n}, Y'_{T_n, T_n, T_n}, Y''_{T_n, T_n, T_n}$
 and $Z_{T_n, T_n, T_n}$ constructed above, respectively.
 Then the Seiberg-Witten invariants of $Z_{n}$ can be easily
 computed using the same technique appeared in the proof of Theorem
 3.3 in~\cite{ss3}.

\begin{theorem}
\label{thm-main-sw}
 The smooth $4$-manifold $Z_n$ has only one $($up to sign$)$
 SW-basic class ${\overline L_n}$ with
 $SW_{Z_n}({\overline L_n}) = \pm n^3$.
\end{theorem}

\begin{proof}
 Since the Seiberg-Witten function  of the $K3$
 surface is $SW_{K3} = 1$ and the Alexander polynomial of $T_n$ is
 $\Delta_{T_n} = nt-(2n-1)+nt^{-1}$, by applying Theorem 1.1
 in~\cite{fs2}, the Seiberg-Witten function of the $V_n=V_{T_n,
 T_n, T_n}$ is equal to
 \[ SW_{V_n} = (n e^{2T} -(2n-1)  + ne^{-2T})^3, \]
 where $T$ is a regular elliptic fiber of $V_n$.
 And the blow-up formula implies that the
 Seiberg-Witten function  of the $V_{n}\sharp 24{\overline{{\mathbf CP}}^2}$
 is also equal to
 \[ SW_{V_n\sharp 24{\overline{{\mathbf CP}}^2}} = SW_{V_n}\cdot
    (e^{E_1} + e^{-E_1})\cdots (e^{E_{24}} + e^{-E_{24}}), \]
 where $E_i$ is an exceptional curve coming from the $i^{th}$ blowing
 up of $V_n$.
 Now applying Theorem~\ref{thm-sw} and Corollary~\ref{cor-sw}
 above for a rational blow-down $4$-manifold $Y'_n$,
 we know that only the following SW-basic classes
 of $V_{n}\sharp 24{\overline{{\mathbf CP}}^2}$
 \[  L_n: = \pm (6T + E_1 + \cdots + E_6 \pm E_7 + E_8 + \cdots + E_{24})\]
 induce SW-basic classes $L'_n$ of $Y'_n$ and they all have
 SW-invariants $SW_{Y'_n}(L'_n) = \pm n^3$.
 Furthermore, applying Theorem 3.2 in~\cite{p1} for $Y'_n$,
 we conclude that such $L'_n$'s are the only SW-basic classes of $Y'_n$.
 Next, again applying Theorem~\ref{thm-sw} and Corollary~\ref{cor-sw}
 above and Theorem 3.2 in~\cite{p1} for the rational blow-down
 $4$-manifold $Z_n$, we conclude that $Z_n$ has only one (up to sign)
 SW-basic class ${\overline L_n}$ which is induced from the following
 SW-basic class of $V_{n}\sharp 24{\overline{{\mathbf CP}}^2}$
 \[  L_n = \pm (6T + E_1 + \cdots + E_6 + E_7 + E_8 + \cdots + E_{24}) \]
 and it has SW-invariant $SW_{Z_n}({\overline L_n}) = \pm n^3$.
\end{proof}

\begin{corollary}
\label{cor-main}
 The $4$-manifolds $Z_n$ $(n \geq 1)$ are simply connected irreducible
 smooth $4$-manifolds which are mutually nondiffeomorphic.
\end{corollary}

\begin{proof}[Proof of Theorem~\ref{thm-main}]
 By Theorem~\ref{thm-main-sc} and Corollary~\ref{cor-main} above,
 the $4$-manifolds $Z_n$ with $n \geq 1$ provide an infinite
 family of smooth $4$-manifolds which are homeomorphic,
 but not diffeomorphic, to
 $3{\mathbf CP}^2 \sharp 8{\overline{{\mathbf CP}}^2}$.
 Hence we are done.
\end{proof}

\bigskip

\section{More Examples}
\label{sec-3}

 In this section we construct more examples of simply connected
 smooth $4$-manifolds with $b_2^+=3$ and $b_2^-=8$ which are
 mutually nondiffeomorphic by using a different configuration.

 First recall that the $4$-manifold $V_{K_1, K_2, K_3}$
 contains a `pseudo-section' with three positive double
 points and with square $-2$.
 At this time, by blowing up at the three positive double points
 of the pseudo-section, we get a $4$-manifold
 $V_{K_1, K_2, K_3}\sharp 3{\overline{{\mathbf CP}}^2}$
 which contains an embedded $2$-sphere of square $-14$
 and which intersects two remaining fishtail fibers
 at the points $p_1$ and $p_2$ respectively.
 Applying $2$ infinitely close blow-ups at $p_1$ and $p_2$
 respectively and blowing up at two positive double points lying
 in the two remaining fishtail fibers, we get a $4$-manifold
 $V_{K_1, K_2, K_3}\sharp 9{\overline{{\mathbf CP}}^2}$
 which contains an embedded $2$-sphere $S$ of square $-18$
 and two copies, say $C_{3,1}$ and $C'_{3,1}$, of a configuration
 $C_{3,1}$ disjoint from $S$.
 And then,
 applying $17$ infinitely close blow-ups at the point $q$,
 we obtain a smooth $4$-manifold
 $V_{K_1, K_2, K_3}\sharp 26{\overline{{\mathbf CP}}^2}$
 which contains three disjoint configurations $C_{3,1}$, $C'_{3,1}$
 and $C_{305, 17}$.
 Since three configurations are disjoint
 in $V_{K_1, K_2, K_3}\sharp 26{\overline{{\mathbf CP}}^2}$,
 we get a new family of smooth $4$-manifolds $\tilde{Y}_{K_1, K_2, K_3}$,
 $\tilde{Y}'_{K_1, K_2, K_3}$ and $\tilde{Z}_{K_1, K_2, K_3}$
 by rationally blowing down along $C_{3,1}$, $C'_{3,1}$ and $C_{305,17}$
 sequentially. Alternatively, we can get a family of smooth $4$-manifolds
 $\tilde{Y}''_{K_1, K_2, K_3}$, $\tilde{Y}'''_{K_1, K_2, K_3}$ and
 $\tilde{Z}_{K_1, K_2, K_3}$ by rationally blowing down along
 $C_{305,17}$, $C_{3,1}$ and $C'_{3,1}$ sequentially.
 Then, using the same techniques as in the proof of
 Theorem~\ref{thm-main-sc} and Theorem~\ref{thm-main-sw} in Section 2,
 we can prove that

\begin{theorem}
\label{thm-extra-sc}
 For any twisted knots $K_1, K_2$ and $K_3$, the $4$-manifold
 $\tilde{Z}_{K_1, K_2, K_3}$ is homeomorphic to
 $3{\mathbf CP}^2 \sharp 8{\overline{{\mathbf CP}}^2}$.
\end{theorem}

\begin{theorem}
\label{thm-extra-sw}
 The smooth $4$-manifold $\tilde{Z}_n = \tilde{Z}_{T_n, T_n, T_n}$
 has only one $($up to sign$)$ SW-basic class ${\overline L_n}$ with
 $SW_{\tilde{Z}_n}({\overline L_n}) = \pm n^3$.
\end{theorem}

\begin{corollary}
\label{cor-extra}
 The $4$-manifolds $\tilde{Z}_n$ $(n \geq 1)$ are simply connected
 irreducible smooth $4$-manifolds with $b_2^+=3$ and $b_2^-=8$
 which are mutually nondiffeomorphic.
\end{corollary}

\noindent {\em Remark.} It is an interesting question to ask
 whether $Z_{K_1, K_2, K_3}$ and $\tilde{Z}_{K_1, K_2, K_3}$
 constructed in this paper are diffeomorphic to the Stipsicz and
 Szab\'{o}'s exotic $4$-manifolds with $b_2^+=3$ and $b_2^-=8$
 in~\cite{ss3}.
 Note that
 we still do not know whether the Stipsicz and Szab\'{o}'s exotic
 $4$-manifolds are simply connected or not. \\


\bigskip

\bigskip



\begin{thebibliography}{BDF}

\bibitem[1]{ch}  A. Casson and J. Harer, {\it Some homology lens spaces which
                 bound rational homology balls}, Pacific. J. Math.
                 {\bf 96} (1981), 23--36
\bibitem[2]{d}   S. Donaldson, {\it Irrationality and the $h$-cobordism conjecture},
                 J. Diff. Geom. {\bf 26} (1987), 141--168
\bibitem[3]{fs1} R. Fintushel and R. Stern, {\it Rational blowdowns of smooth
                 $4$-manifolds}, J. Diff. Geom. {\bf 46} (1997), 181--235
\bibitem[4]{fs2} R. Fintushel and R. Stern, {\it Knots, links and
                 $4$-manifolds}, Invent. Math. {\bf 134} (1998), 363--400
\bibitem[5]{fs3} R. Fintushel and R. Stern, {\it Double node neighborhoods and
                 families of simply connected $4$-manifolds with $b^+=1$},
                 to appear in J. Amer. Math. Soc. (arXiv:math.GT/0412126)
\bibitem[6]{fm}  R. Friedman and J. Morgan, {\it Smooth $4$-manifolds and complex
                 surfaces}, Ergebnisseder Mathematik und ihrer Grenzgebiete
                 {\bf 27}, Springer-Verlag, 1994
\bibitem [7]{g}  R. Gompf, {\it A new construction of symplectic  manifolds},
                 Annals of Math. {\bf 142} (1995),  527--595
\bibitem[8]{hkk} J. Harer, A. Kas and R. Kirby, {\it Handlebody decompositions
                 of complex surfaces}, Memoirs of Amer. Math. Soc. {\bf 62}, 1986
\bibitem[9]{k}   D. Kotschick, {\it On manifolds homeomorphic to
                 ${\mathbf CP}^2 \sharp 8{\overline{{\mathbf CP}}^2}$},
                 Invent. Math. {\bf 95} (1989), 591--600
\bibitem[10]{mo} J. Morgan and K. O'Grady, {\it Differential topology of complex
                 surfaces. Elliptic surfaces with $p_g =1$: smooth classification},
                 Lecture Notes in Mathematics {\bf 1545}, Springer-Verlag, 1993
\bibitem[11]{pd1} D. Park, {\it Exotic smooth structures on
                  ${\mathbf CP}^2 \sharp n{\overline{{\mathbf CP}}^2}$},
                  Proc. Amer. math. Soc. {\bf 128} (2000), 3057--3065
\bibitem[12]{pd2} D. Park, {\it Exotic smooth structures on
                  ${\mathbf CP}^2 \sharp n{\overline{{\mathbf CP}}^2}$, Part II},
                  Proc. Amer. math. Soc. {\bf 128} (2000), 3067--3073
\bibitem[13]{pd3} D. Park, {\it Constructing infinitely many smooth structures
                  on $3{\mathbf CP}^2 \sharp n{\overline{{\mathbf CP}}^2}$},
                  Math. Ann. {\bf 322} (2000), 267--278
\bibitem[14]{p1}  J. Park, {\it Seiberg--Witten invariants of generalized rational
                  blow-downs}, Bull. Austral. Math. Soc. {\bf 56} (1997),
                  363--384
\bibitem[15]{p2}  J. Park, {\it Exotic smooth structures on 4-manifolds},
                  Forum Math. {\bf 14} (2002), 915--929
\bibitem[16]{p3}  J. Park, {\it Exotic smooth structures on 4-manifolds, II},
                  Topology Appl. {\bf 132} (2003), 195--202
\bibitem[17]{p4}  J. Park, {\it Simply connected symplectic 4-manifolds with $b_2^+=1$
                  and $c_1^2=2$}, Invent. Math. {\bf 159} (2005), 657--667
\bibitem[18]{pss} J. Park, A. Stipsicz and Z. Szab\'{o}, {\it Exotic smooth structures
                  on ${\mathbf CP}^2 \sharp 5{\overline{{\mathbf CP}}^2}$},
                  Math. Res. Letters {\bf 12} (2005), 701--712
\bibitem[19]{st}  A. Stipsicz, {\it Donaldson series and ($-1$)-tori},
                  J. Reine Angew. Math. {\bf 465} (1995), 133--144
\bibitem[20]{ss1} A. Stipsicz and Z. Szab\'{o}, {\it The smooth classification
                  of elliptic surfaces with $b_2^+ >1$},
                  Duke Math. J. {\bf 75} (1994), 1--50
\bibitem[21]{ss2} A. Stipsicz and Z. Szab\'{o}, {\it An exotic smooth structure
                  on ${\mathbf CP}^2 \sharp 6{\overline{{\mathbf CP}}^2}$},
                  Geometry \& Topology {\bf 9} (2005),  813--832
\bibitem[22]{ss3} A. Stipsicz and Z. Szab\'{o}, {\it Small exotic $4$-manifolds
                  with $b_2^+=3$}, arXiv:math.GT/0501273
\bibitem[23]{sz}  Z. Szab\'{o}, {\it Irreducible four-manifolds with small Euler
                  characteristics}, Topology {\bf 35} (1996), 411--426

\end{thebibliography}
\end{document}